\documentclass[12pt]{amsart}
\usepackage[T1]{fontenc}
\usepackage[utf8]{inputenc}
\usepackage{fourier}
\usepackage[active]{srcltx}
\usepackage{a4wide}
\usepackage{amsthm,amsfonts,amsmath,mathrsfs,amssymb}
\usepackage{dsfont}
\def\be{\begin{equation}}
\def\ee{\end{equation}}

\newtheorem*{property*}{Property}

\newtheorem{conjecture}{Conjecture}
\newtheorem*{conjecture*}{Conjecture}
\newtheorem{remark}{Remark}
\newtheorem*{completeness*}{Completeness property}
\newtheorem*{theorem*}{Theorem}
\newtheorem{theorem}{Theorem}
\newtheorem{proposition}{Proposition}
\newtheorem*{proposition*}{Proposition}
\newtheorem{lemma}{Lemma}
\newtheorem{corollary}{Corollary}
\theoremstyle{remark}

\newcommand{\B}{{\mathbf{B}}}
\newcommand{\nc}{\newcommand}

\newcommand{\cc}{{\mathbf Y}}
\newcommand{\cl}{{\mathbf Y_{\ell}}}

\newcommand{\PP}{\mathscr{P}}

\newcommand{\RR}{R}

\newcommand{\M}{\mathcal{M}}

\newcommand{\N}{{\mathbb N}}

\nc{\supp}{\operatorname{supp}}
\nc{\Real}{\operatorname{Re}}
\nc{\Imag}{\operatorname{Im}}
\nc{\dif}{\operatorname{d}} \nc{\im}{\operatorname{i}}

\nc{\Hi}{{\mathscr{H}}^\infty} \nc{\Ht}{{\mathscr{H}}^2}
\nc{\Hone}{{\mathscr{H}}^1} \nc{\ol}{\overline} \nc{\bz}{\mathbf{z}}
\nc{\bw}{\mathbf{w}} \nc{\eps}{\varepsilon}

\allowdisplaybreaks[3]

\begin{document}

\title[Extreme Values of Derivatives of zeta and $L$-functions]
{Extreme Values of Derivatives   of zeta and $L$-functions}
\author{Daodao Yang}
\address{Institute of Analysis and Number Theory \\ Graz University of Technology \\ Kopernikusgasse 24/II, 
A-8010 Graz \\ Austria}

\email{yang@tugraz.at \quad yangdao2@126.com}

\dedicatory{Dedicated to   Kristian Seip
  on the occasion of his 60th birthday}

\begin{abstract}
It is proved  that  as $T \to \infty$,  uniformly for all positive integers $\ell \leqslant (\log_3 T) / (\log_4 T)$, we have
\begin{equation*}  
\max_{T\leqslant t\leqslant 2T}\left|\zeta^{(\ell)}\Big(1+it\Big)\right| \geqslant \big(\cl + o\left(1\right)\big)\left(\log_2 T \right)^{\ell+1} \,,
\end{equation*}
where \,$\cl = \int_0^{\infty} u^{\ell} \rho (u) du$. Here $\rho(u)$ is the Dickman function. We have \,$\cl > e^{\gamma}/(\ell + 1)$ and \,$ \log\, \cl = \left(1 + o\left(1\right) \right) \ell \log \ell$ when $ \ell \to \infty $, which significantly improves
previous results in \cite{DongNote, DY}. Similar results are established for Dirichlet $L$-functions. On the other hand, when assuming the Riemann Hypothesis and the Generalized Riemann Hypothesis, we establish upper bounds for $	\left| \zeta^{(\ell)}\left(1+it\right)\right| $ and $\left|L^{(\ell)}(1, \chi) \right|$.  Furthermore,  when assuming the   Granville-Soundararajan  Conjecture is true, we establish the following asymptotic formulas 
$$\max_{ \substack{  \chi \neq \chi_0 \\ \chi(\text{mod}\, q)}} \left|L^{(\ell)}(1, \chi) \right| \sim \cl\left(\log_2 q\right)^{\ell+1},\,\, \quad \text{as}\,\quad q \to \infty,$$
where  $q$  is prime and $\ell \in \mathbb{N}$ is given.
\end{abstract}
\maketitle
\section{Introduction}
This paper establishes the following  results for extreme values of derivatives of the Riemann
zeta function on the 1-line.   Throughout the paper, we define $\cl = \int_0^{\infty} u^{\ell} \rho (u) du$  and $\rho(u)$ denotes  the Dickman function.

\begin{theorem}\label{Main: sigma =  1} 
As $T \to \infty$,  uniformly for all positive integers $\ell \leqslant (\log_3 T) / (\log_4 T)$, we have
\[  \max_{T\leqslant t\leqslant 2T}\left|\zeta^{(\ell)}\Big(1+it\Big)\right| \geqslant \big(\cl + o\left(1\right)\big)\left(\log_2 T \right)^{\ell+1}  \,. \]

\end{theorem}

\begin{remark}\label{compute}
By taking derivatives for  the Laplace transform (see Lemma\,\,\ref{Laplace}) of the Dickman function and applying Fa\`a di Bruno's formula (for example, see \cite[page 134--137]{Comb}), we can obtain a formula for $\cl$ in terms of Bell polynomials. 
Namely, let \,$\B_{\ell}\left(x_1, x_2, \cdots, x_{\ell}\right)$ be the  $\ell$-th complete exponential Bell polynomial, then we have \,$\cl = e^{\gamma} (-1)^{\ell} \B_{\ell}\left( -1, \frac{1}{2}, \cdots, \frac{(-1)^{\ell}}{\ell}\right)$\,. For instance, from $\B_1(x_1) = x_1$,  $\B_2(x_1, x_2) = x_1^2 + x_2$, $\B_3(x_1, x_2, x_3) = x_1^3 + 3 x_1 x_2 + x_3$, we can compute \,\,$\cc_{1} = e^{\gamma}$,  $\cc_{2} = 3e^{\gamma}/2$ and $\cc_{3} = 17e^{\gamma}/6$.
\end{remark}

\begin{remark}
By the asymptotic formula \eqref{decay}, we have\, $ \log\, \cl = \left(1 + o\left(1\right) \right) \ell \log \ell$, as \,$ \ell \to \infty $.
\end{remark}

Our result on the Riemann zeta function can  be generalized to $L$-functions. In the following theorem, we consider  the case of Dirichlet $L$-function $L(s, \chi)$ associated with  non-principal characters $\chi$(mod $q$).

\begin{theorem}\label{Lchi_lower}
	Let  $q$ be prime,  then uniformly for all positive integers  $\ell \leqslant \log_3 q / \log_4 q$, we have
\begin{align}\label{ExtremeL}
   \max_{ \substack{  \chi \neq \chi_0 \\ \chi(\text{mod}\, q)}} \left|L^{(\ell)}(1, \chi) \right| \geqslant \left(\cl+o\left(1\right)\right) \left(\log_2 q\right)^{\ell+1},\,\, \quad \text{as}\,\quad q \to \infty.
\end{align}  
\end{theorem}

\begin{remark}
1) The above result does not hold for general moduli $q$. For instance, assume $q = \left(\prod_{p \leqslant X} p \right) \cdot m $\,, with $m \in \N$ and $X = \frac{1}{2} \log q$. This assumption will force $\chi(k) = 0$ if $\,\,\exists p \leqslant X,  p|k $ and thus will make $\left|L^{(\ell)}(1, \chi) \right|$ small.  2) However, if $q$ is not divisible by small primes (for instance, consider the case that any prime factor of $q$ is larger than $q^{\frac{1}{10}}$), then the above theorem will still hold. For simplicity, we state the result for prime moduli.
\end{remark}
For upper bounds of $	\left| \zeta^{(\ell)}\left(1+it\right)\right| $ and $\left|L^{(\ell)}(1, \chi) \right|$, we have following two results when assuming the Riemann Hypothesis (RH) and the Generalized Riemann Hypothesis (GRH). 

\begin{theorem}\label{RH_upperBound}
Fix $\ell \in \mathbb{N}$. Assuming  RH, we have 
	\begin{align*}
		\left| \zeta^{(\ell)}\left(1+it\right)\right| \leqslant \left( 2^{\ell + 1}\, \cl + o(1) \right) \left(\log_2 t \right)^{\ell + 1}\,, \quad \text{as}\quad t \to \infty\,.
	\end{align*}
\end{theorem}

\begin{theorem}\label{GRH_upperBound}
Fix $\ell \in \mathbb{N}$. Let $\chi$ be any non-principal character (mod $q$), and assume   GRH for $L(s, \chi)$. Then
	\begin{align*}
	\left|L^{(\ell)}(1, \chi) \right| \leqslant \left( 2^{\ell + 1}\, \cl + o(1) \right) \left(\log_2 q \right)^{\ell + 1}\,, \quad \text{as}\quad q \to \infty\,.
	\end{align*}
\end{theorem}

The key ingredient to prove Theorem \ref{GRH_upperBound} is the following theorem of   Granville and Soundararajan \cite[Theorem 2]{largeGS}. 

To state their result, we need some definitions.  Let $f$ be an arithmetic function. Define the functions $\Psi(x,y)$ and $ \Psi(x,y; f)$ as
\begin{align*}
 \Psi(x,y):\, = \sum_{\substack{n \leqslant x \\ P^{+}(n) \leqslant y}}   1\,, \quad \quad \quad \Psi(x,y; f):\, = \sum_{\substack{n \leqslant x \\ P^{+}(n) \leqslant y}}   f(n)\,.
\end{align*}

\begin{theorem*}[Granville-Soundararajan]\label{GRH_GS}
Let $\chi$ be any non-principal character (mod $q$), and  assume  the Riemann Hypothesis for $L(s, \chi)$.
If $1\leqslant x \leqslant q$ and $y\geqslant \log^2q\log^{2} x (\log_2 q )^{12}$,  then
$$
\sum_{n\leqslant x} \chi(n) = \Psi(x,y;\chi)+O\biggl( \frac{\Psi(x,y)}
{(\log_2 q )^2}\biggr)\,.
$$
Further
$$
\biggl| \sum_{n\leqslant x} \chi(n) \biggr| \ll
\Psi(x, \log^{2}q (\log_2 q )^{20}),
$$
and so the following estimate  holds 
\begin{align*}
\biggl |\sum_{n\leqslant x} \chi(n)\biggr| = o(x) \,,   
\end{align*}
when $\log x/\log_2 q  \to \infty$ as $q\to \infty$.

\end{theorem*}
When $g(n) = n^{-it}, \,\forall n \in \N$, we  write  $\Psi(x,y; t)$ in place of $\Psi(x,y; g)$. 
Then we have the following result analogous to the Granville-Soundararajan Theorem.  
\begin{theorem}\label{RH_DY}
Assume  RH and let $T$ be sufficiently large.
If \,$2\leqslant x \leqslant T$, $ T + y + 3\leqslant t \leqslant  T^{1000}$  and 
$ y\geqslant \log^2 T\log^{2} x (\log_2 T)^{12}$,  then
\begin{align}\label{Approx_BySmooth}
\sum_{n\leqslant x}  \frac{1}{n^{it}} = \Psi(x,y;t)+O\biggl( \frac{\Psi(x,y)}
{(\log_2 T)^2}\biggr)\,.
\end{align}
Further, 
\begin{align}\label{Upper_BySmooth}
    \biggl| \sum_{n\leqslant x}  \frac{1}{n^{it}}  \biggr| \ll
\Psi(x, \log^{2}T (\log_2 T)^{20}),\quad \forall x \in [2,\, T]\,, \forall t \in [T + \log^{2}T (\log_2 T)^{15} ,\,  T^{1000}]\,, 
\end{align}
and so the following estimate  holds 
\begin{align*}
\biggl |\sum_{n\leqslant x} \frac{1}{n^{it}}\biggr| = o(x) \,,   \quad \forall x \in [2,\, T]\,, \forall t \in [T + \log^{2}T (\log_2 T)^{15} ,\,  T^{1000}]\,, 
\end{align*}
when $\log x/\log_2 T \to \infty$ as $T\to \infty$.

\end{theorem}

In \cite{largeGS}, Granville and Soundararajan also made the following
conjecture.

\begin{conjecture*}[Granville-Soundararajan]\label{GSC}
There exists a constant $A>0$ such that
for any non-principal character $\chi$ (mod $q$), and for any
$1\leqslant x\leqslant q$ we have, uniformly,  
$$
\sum_{n\leqslant x} \chi(n) = \Psi(x,y;\chi)+o(\Psi(x,y;\chi_0)),
$$
where $y= (\log q +\log^2 x) (\log_2 q )^A$.
\end{conjecture*}

A consequence of  the   Granville-Soundararajan  Conjecture  is that the constant appearing in Theorem \ref{Lchi_lower} is sharp. 

\begin{theorem}\label{GSC_upper}
Assume  the   Granville-Soundararajan  Conjecture is true. Fix $\ell \in \mathbb{N}$. Let $\chi$ be any non-principal character (mod $q$),  then
	\begin{align*}
	\left|L^{(\ell)}(1, \chi) \right| \leqslant \left(  \cl + o(1) \right) \left(\log_2 q \right)^{\ell + 1}\,, \quad \text{as}\quad q \to \infty\,.
	\end{align*}
\end{theorem}
Combining with the lower bound, we immediately get the following 
asymptotic formulas.
\begin{corollary}\label{GSC_Cor}
Assume  the   Granville-Soundararajan  Conjecture is true. Fix $\ell \in \mathbb{N}$. Let  $q$ be prime, then
	\begin{align*}
	\max_{ \substack{  \chi \neq \chi_0 \\ \chi(\text{mod}\, q)}} \left|L^{(\ell)}(1, \chi) \right| \sim \cl\left(\log_2 q\right)^{\ell+1},\,\, \quad \text{as}\,\quad q \to \infty.
	\end{align*}
\end{corollary}

We have the following analogous conjecture, which can imply that the constant appearing in Theorem \ref{Main: sigma =  1} is sharp. 

\begin{conjecture}\label{DYC}
There exists a constant $A>0$ such that
 for any
$1\leqslant x \leqslant T$, $2 T \leqslant t \leqslant 5T$, we have, uniformly,  
$$
\sum_{n\leqslant x}  \frac{1}{n^{it}} = \Psi(x,y;t)+o\left( \Psi(x,y)
\right)\,,  \quad \text{as}\quad T \to \infty\,,
$$
where $y= (\log T +\log^2 x) (\log_2 T)^A$.
\end{conjecture}

\begin{theorem}\label{DYC_Cor}
Assume  Conjecture \ref{DYC} is true. Fix $\ell \in \mathbb{N}$. Then
	\[  \max_{T\leqslant t\leqslant 2T}\left|\zeta^{(\ell)}\Big(1+it\Big)\right| \sim \cl \left(\log_2 T \right)^{\ell+1}  ,  \quad \text{as}\quad T \to \infty\,. \]
\end{theorem}

The problem of obtaining extreme values of $\left|\zeta\left(1+it\right)\right|$ was first considered by Bohr and Landau, who  established the result $\zeta(1+it) = \Omega(\log_2 t)$ (see \cite[Thm 8.5]{T}) in 1910.  In 1924,  Littlewood  (see \cite[Thm 8.9(A)]{T})  improved the result of  Bohr and Landau, by proving  that $\limsup_{t\to\infty} |\zeta(1+it)|/(\log_2 t) \geqslant e^{\gamma}.$  Littlewood's result has been improved in the past century by Levinson \cite{L},  by 
Granville-Soundararajan \cite{GS}, and by Aistleitner-Mahatab-Munsch\cite{AMM},
who established that 
$\max_{\sqrt T \leqslant t\leqslant T} \left|\zeta(1 + it)\right|  \geqslant e^{\gamma} (\log_2 T + \log_3 T + C),$
for  some constant $C$  which can be effectively computed.  
Littlewood also established conditional results for the upper bound of $|\zeta(1+it)|$. When assuming the truth of the Riemann hypothesis (RH),  he proved that $|\zeta(1+it)| \leqslant (2 e^{\gamma} + o(1)) \log_2 t\,, $ as $t\to\infty$  (see \cite[Thm 14.9]{T}). Furthermore, he conjectured that the maximum of $|\zeta(1+it)|$ on the interval  $[1, T]$ should satisfy the asymptotic formula $ \max_{1 \leqslant t\leqslant T}  \left|\zeta(1 + it)\right| \sim e^{\gamma}\log_2 T .$   In \cite{GS},  Granville and  Soundararajan  made the  stronger conjecture that\,
$
  \max_{T \leqslant t\leqslant 2T}  \left|\zeta(1 + it)\right| = e^{\gamma}(\log_2 T + \log_3 T + C) + o(1), 
$
for some constant $C$ which is also effectively computable.

When  $\sigma \in [1/2, \,1)$, the problem of  obtaining  extreme values of  $\left|\zeta(\sigma+it)\right|$  also has a long history. For more background and results, see the recent survey \cite{Sound} and \cite{A, BR, BS1, BS2, BSNote, CC, CSo, delaBT, Dong, GS, Hi, M, So, V}. Here we mention the recent breakthrough result by Bondarenko and  Seip \cite{BS1, BS2} who prove that:\[ \max_{1 \leqslant t \leqslant T} \left|\zeta\Big(\frac{1}{2}+it\Big)\right| \geqslant \exp\left(c\sqrt{\frac{\log T \log_3 T}{\log_2 T}}\right), \quad \forall T \gg 1\,,\] 
for any constant $c < 1$.
After refining methods of  Bondarenko-Seip, de la Bret\`eche and Tenenbaum \cite{delaBT}  show that  any $c < \sqrt 2$ is permissible in the above result.

In the past five years, the problem of obtaining extreme values of  derivatives of the Riemann zeta function have been studied. In \cite{Kalmynin}, Kalmynin obtained $\Omega$-results for the Riemann zeta function and its derivatives $\zeta^{(\ell)}(\sigma+it)$, when $\sigma = \sigma(t) \to 1^{-}, $ as $t\to \infty$. 

In \cite{DY}, we established $\Omega$-results for $\big|\zeta^{(\ell)}(\sigma+it)\big|$ when $\ell \in \mathbb N $
 and $\sigma \in [1/2, \,1)$ are  given.   These results are comparable with the best currently known lower bounds for maximum of $\big|\zeta(\sigma+it)\big|$.   When $\sigma = 1$,  we obtain  lower bounds different from the case of $\big|\zeta(1+it)\big|$. Namely,  we  established that
$ 
\max_{T\leqslant t\leqslant 2T}\left|\zeta^{(\ell)}\left(1+it\right)\right| \geqslant e^{\gamma}\cdot \ell^{\ell}\cdot (\ell+1)^{ -(\ell+1)}\cdot\left(\log_2 T - \log_3 T + O(1)\right)^{\ell+1} \,,
$ uniformly for all positive integers $\ell \leqslant (\log T) / (\log_2 T)$, when $T$ is sufficiently large. On the other hand,  in \cite{Daodao} we proved that on RH,  $~\vert \zeta^{(\ell)}\left(1+it\right) \vert \ll_{\ell} (\log_2 t)^{\ell+1}~$
for sufficiently large $t$, where the implied constants are effectively computable.  Refining methods of \cite{DY},  Z. Dong and B. Wei \cite{DongNote} proved that $ 
\max_{T\leqslant t\leqslant 2T}\left|\zeta^{(\ell)}\left(1+it\right)\right| \geqslant \left(e^{\gamma}/(\ell+1) + o(1) \right)\cdot\left(\log_2 T \right)^{\ell+1} \,,
$ uniformly for all positive integers $\ell \leqslant (\log T) / (\log_2 T)$, as $T \to \infty$. The constant $e^{\gamma}/(\ell+1)$ improve the constant $e^{\gamma}\cdot \ell^{\ell}\cdot (\ell+1)^{ -(\ell+1)}$ by a factor $(1 + 1/\ell)^{\ell}$. However, we still have $\lim_{\ell \to \infty} e^{\gamma}/(\ell+1) = 0$. In contrast, in our new result,  we have \,$\lim_{\ell \to \infty}\cl = \infty $.  Also we have \,$\cl > e^{\gamma}/(\ell + 1), \forall \ell > 0 $. This is due to the following  identity $$\cl - \frac{e^{\gamma}}{\ell + 1} = \int_1^{\infty} \left(u^{\ell}-\frac{1}{\ell + 1}\right) \rho(u) du \,,$$
and the fact that $\rho(u)$ is always positive. 

By Theorem \ref{RH_upperBound},  assuming  RH, we have  
$\left|\zeta^{\prime \prime}\left(1+it\right)\right| \leqslant \left( 12 e^{\gamma} + o(1) \right) \left(\log_2 t\right)^{3} $  and $\left|\zeta^{(3)}\left(1+it\right)\right| \leqslant \left( \frac{136}{3} e^{\gamma} + o(1) \right)  \left(\log_2 t\right)^{4} $, which improve corresponding results in \cite{Daodao}.

The study of extreme  values of $L$-functions is  an important problem in analytic number theory.  Given a negative fundamental discriminant $d$, one can associate a primitive Dirichlet character $\chi_d$ (mod $|d|$) by defining $\chi_d(n) = \left(\frac{d}{n}\right)$, using the Kronecker–Legendre symbol.  The value $L(1, \chi_d)$ is related to the class number of $\mathbb Q(\sqrt{d})$ via  Dirichlet’s
class number formula\begin{align*}
   L(1, \chi_d) = \frac{2\pi h }{\omega \sqrt{|d|}} \,,
\end{align*}where $h$ is the class number of $\mathbb Q(\sqrt{d})$, and $\omega$ 
denotes the number of roots of unity in $\mathbb Q(\sqrt{d})$.

Let $\chi$ be any non-principal character (mod $q$). Assuming GRH,  Littlewood \cite{Littlewood} 
proved that 
$$
L(1,\chi) \sim \prod_{p\leqslant \log^2 q} \Big(1-\frac{\chi(p)}{p}\Big)^{-1}, \quad \text{as}\quad q \to \infty\,,
 $$
from which one can immediatley obtain that
$
\left|L(1,\chi) \right|\leqslant \left(2 e^{\gamma}+o(1)\right) \log_2 q 
$ by Mertens'  theorem. In another direction, Chowla \cite{Chowla} showed that 
there exist arbitrarily large $q$ and non-principal  characters $\chi$(mod $q$)
such that 
$
\left|L(1,\chi) \right| \geqslant \left(e^{\gamma}+o(1)\right) \log_2 q 
$.   As in the proof of Theorem \ref{GSC_upper},  the   Granville-Soundararajan Conjecture   implies $
\left|L(1,\chi) \right|\leqslant \left( e^{\gamma}+o(1)\right) \log_2 q 
$ for non-principal  characters $\chi$.
We mention three best  results known for  $L(1, \chi)$.  For further information and results about extreme  values of $L$-functions, we refer to the survey \cite{Sound} and \cite{AKMP, delaBT, GSchi, GSchid, GS, LL, So}.
In \cite{GS},  Granville and Soundararajan  established
that for sufficiently large prime $q$ and any given $A \geqslant 10$ there are at least 
$q^{1-1/A}$ characters
$\chi$(mod $q$) for which
$$\left|L(1, \chi)   \right| \geqslant e^{\gamma}  \left( \log_2 q + \log_3 q - \log_4 q  - \log A - C  \right)\,,$$
for some absolute constant $C$. In \cite{AKMP}, Aistleitner,  Mahatab,  Munsch and Peyrot  proved that when fix $\epsilon > 0$, then for all sufficiently large prime $q$, we have $$\max_{\chi \neq \chi_0}\left|L(1, \chi)   \right| \geqslant e^{\gamma}  \left( \log_2 q + \log_3 q - (1+\log_2 4) -\epsilon  \right)\,.
$$
In \cite{LL}, when assuming GRH, Lamzouri-X. Li-Soundararajan obtained the following upper bound for  primitive
character $\chi $ modulo $q$
\begin{align*}
	\left|L(1, \chi)   \right| \leqslant 2 e^{\gamma}  \left(\log_2 q  - \log 2  + \frac{1}{2} + \frac{1}{\log_2 q }\right)\,, \quad \forall q \geqslant 10^{10}\,.
	\end{align*}

Like  $L(1, \chi_d)$, the value of the first derivative  $L^{\prime}(1, \chi_d)$ can be related to the class number as well, namely, via the following Chowla-Selberg formula \cite[page 110]{SC}
\begin{align*}
  L^{\prime}(1, \chi_d) =   -\frac{\pi}{\sqrt{|d|}}\sum_{m = 1}^{|d|}\chi_d(m) \log \Gamma \left(  \frac{m}{|d|}\right) + \frac{2h\pi(\gamma + \log 2\pi)}{\omega \sqrt{|d|}}\,.
\end{align*}
In \cite[page 524]{IK}, Iwaniec and  Kowalski mention that when assuming  GRH for $L(s, \chi_d)$, one can obtain that $\left|L^{(\ell)}(1, \chi_d) \right|\ll (\log_2 |d|)^{\ell + 1}$. However, they do not point out what the implicit constants could be. On the other hand, we do not find literatures on large values of $\left|L^{(\ell)}(1, \chi) \right|$. Theorem \ref{Lchi_lower} and \ref{GRH_upperBound} can be considered as generalizations of theorems  of Littlewood  and Chowla.

The study of character sums is another central problem in number theory.  In many cases, one would like to know when the following character sum is $o(x)$, 
$$\sum_{n \leqslant x} \chi(n)\,,$$
   where $\chi$ is a non-principal Dirichlet character $\chi$(mod $q$).

In \cite{MV},   Montgomery and Vaughan show that the above character sums can be conditionally approximated by character sums over integers with small prime factors. More precisely, they 
prove that if $\chi$(mod $q$) is non-principal and GRH holds then
$$
\sum_{n\leqslant x} \chi(n) = \Psi(x, y; \chi) + O(xy^{-\frac 12} \log^4 q),
$$
when $\log^4 q \leqslant y\leqslant x\leqslant q$. One of main results in \cite{MV} states that on GRH,  
\begin{align}\label{GRHPV}
     \left|\sum_{n\leqslant x} \chi(n)\right| \ll \sqrt q \log_2 q 
\,,\end{align}
for any non-principal character $\chi$ modulo $q$ and any $x$. On GRH, Granville and Soundararajan \cite{largeGS2007} find an implicit constant in \eqref{GRHPV} for  primitive
character $\chi $ modulo $q$. The upper bound
 \eqref{GRHPV} can be used to  improve the error term in the approximation formula \eqref{approxD} for $L^{(\ell)}(1, \chi)$. In \cite{largeGS}, Granville and Soundararajan   refine the methods of Montgomery-Vaughan  to 
obtain the  result mentioned early in the paper, which turns out to be a key to our understanding of upper bounds of $\left|L^{(\ell)}(1, \chi) \right|$. And the Theorem \ref{RH_DY} are based on  the work of Montgomery-Vaughan and Granville-Soundararajan, in particular following methods of Granville-Soundararajan.

We will  use  Soundararajan's  resonance methods  \cite{So}  to prove  Theorem \ref{Main: sigma =  1} and Theorem \ref{Lchi_lower}. The key ingredient  is the following  Proposition \ref{maxRatio}.

\begin{proposition}\label{maxRatio}
As $T \to \infty$,  uniformly for all positive numbers $\ell \leqslant (\log_3 T) / (\log_4 T)$, we have
\begin{equation*}
\sup_{r} \Big|\sum_{mk = n\leqslant \sqrt T} \frac{r(m)\overline{r(n)}}{ k} (\log k)^{\ell}
\Big| \Big/ \Big(\sum_{n\leqslant \sqrt T }|r(n)|^2\Big) 
\geqslant \big(\cl + o\left(1\right)\big)\left(\log_2 T \right)^{\ell+1} ,
\end{equation*}
where the supremum is taken over all functions $r :\, \mathbb N \to \mathbb C$ satisfying that the denominator is not equal to zero, when the parameter $T$ is given. 
\end{proposition}

\textbf{Notations:} in this paper,  $\gamma$ denotes the Euler constant. We write $\log_j$
for the $j$-th iterated logarithm, so
for example  $\,\log_2 T \,= \log\log T,$  $\log_3 T \,= \log\log\log T$. $P^{+}(n)$ denotes the largest prime factor of $n$. $p$ denotes a prime number and $p_n$ denotes the $n$-th prime.


\section{Preliminary Results}
Recall that the  function $\Psi(x, y) = \#\left\{  n \leqslant x \Big| P^{+}(n) \leqslant y \right\} $ counts the number of integers $n$ not exceed $x$ with prime factors at most $y$.  The Dickman function $\rho(u)$ is a continuous function defined by the initial condition  $ \rho (u)=1 $ for  $ 0 \leqslant u \leqslant 1 $ and satisfies the following differential equation
\begin{align}\label{eqDickman}
 u\rho^{\prime}(u)+\rho (u-1)=0 \,,\,\quad u>1 \,.   
\end{align}
From the definition, the Dickman function $\rho(u)$ is a positive decreasing function.  In 1930, Dickman   \cite{Dickman} proved that for fixed $u > 0$,  $\lim_{x \to \infty}  \Psi(x, x^{\frac{1}{u}})/x$  exists and equals to $\rho(u)$. We will use the following strong form of this asymptotic formula and an asymptotic formula for $\rho(u)$. In the following lemma, \eqref{JianJin} is due to Hildebrand\cite{H1986a}. The upper bound of \eqref{JianJin2} is due to de Bruijn\cite{deB1966}, while the lower bound of \eqref{JianJin2} is due to  Hildebrand\cite{H1986a}. And the asymptotic formula \eqref{decay}  for $\rho(u)$   was obtained by de Bruijn\cite{deB1951}.
\begin{lemma}[Thm 1.1, 1.2,\,Cor 2.3\, \cite{HT}]\label{HTDickman}
Let $x \geqslant y \geqslant 2$ be real numbers, and put $ u = \frac{\log x}{\log y} $. For any fixed $\epsilon > 0$ the asymptotic formula
\begin{align}\label{JianJin}
    \Psi(x , y) = x \rho (u) \left( 1 +  O \left( \frac{\log(u+1)}{\log y}    \right)  \right)
\end{align}
holds uniformly in the range $1 \leqslant u \leqslant \exp\left(  \left(   \log y \right)^{\frac{3}{5} - \epsilon}  \right)$.\, The weaker relation
\begin{align}\label{JianJin2}
\log \frac{\Psi(x,y)}{x} = \left( 1+ O\left(\exp(-(\log u)^{\frac 35 
-\epsilon})\right)\right)\log \rho(u) 
\end{align}
holds uniformly in the range $1\leqslant u\leqslant y^{1-\epsilon}$.  And as $u \to \infty$,
\begin{align}\label{decay}
    \log \rho (u) = - u \left(\log u + \log_2 (u+2) - 1 + O\left( \frac{\log_2 (u+2)}{\log (u+2)} \right)\right)\,.
\end{align}
\end{lemma}
The following lemma is on the  Laplace transform of the Dickman function, which is useful for us to compute $\cl$, as mentioned in Remark \ref{compute}.
\begin{lemma}[Lemma 2.6\,\cite{HT}, \, Thm 7.10\,\cite{M2}]\label{Laplace}
For any real or complex number $s$ we have
\begin{align*}
  \int_{0}^{\infty} \rho(u) e^{-us}du=\exp\left(\gamma + \int_0^{s} \frac{e^{-z}-1}{z}dz\right)\,.
\end{align*}
\end{lemma}

We have the following conditional approximation formula for $\log \zeta(\sigma+it)$, which is adapted from Lemma 1 of \cite{GS}.

\begin{lemma}[Granville-Soundararajan]\label{approx_logZeta}
Assume RH. Let $y\geqslant 2$ and $t \geqslant y+3$. For $\frac{1}{2} < \sigma \leqslant 1$, we have
$$
\log \zeta(\sigma+it)= \sum_{n=2}^{[y]}
\frac{\Lambda(n)}{n^{\sigma+it} \log n} + O\Big(
\frac{\log t}{(\sigma_1-\frac{1}{2})^2}y^{\sigma_1-\sigma}\Big),
$$
where we put $\sigma_1 = \min(\frac{1}{2} +\frac{1}{\log y},
\frac{\sigma}{2} + \frac{1}{4})$.
 
\end{lemma}
We have the following unconditional approximation formula for $\zeta^{(\ell)}(\sigma+ it)$. The constant 6.28 can be replaced by any positive number smaller than $2\pi$ (see \cite[Lemma 2]{HaL} ).

\begin{lemma}[Lemma 1 \cite{DY}]
 Let  $\sigma_0 \in (0, 1)$ be fixed. If $T$ is sufficiently large,
then uniformly for $\epsilon >0$, $t \in [T, \, 6.28T]$, $\sigma \in[ \sigma_0 + \epsilon, \, \infty)$ and all positive integers $\ell$, we have
\begin{equation}\label{approx}
   (-1)^{\ell} \zeta^{(\ell)}(\sigma+ it) = \sum_{n \leqslant T} \frac{(\log n)^{\ell}}{n^{\sigma+it}} + O\Big(  \frac{\ell !}{\epsilon^{\ell}}\cdot T^{-\sigma+\epsilon}\Big)\,,
\end{equation}
where the implied constant in big $O(\cdot)$ only depends on $\sigma_0$ .
\end{lemma}
\section{Proof of  Proposition \ref{maxRatio}}

\begin{proof}

Let $T$ be large. Let $w = \pi(y)$. Define $y$, $b$ and $\PP(y, b)$  as  follows
\[ y = \frac{\log T}{3  (\log_2 T)^3}    \,,\quad b = [(\log_2 T)^3]\,, \quad    \PP(y, b) = \prod_{p \leqslant y} p^{b-1} = \prod_{i = 1}^w p_i^{b-1}\,\;.\]

Note that $\PP(y, b) \leqslant \sqrt T$. Let $\M$ be the set of divisors of $\PP(y, b)$ and define the function $r:\,\N \to \{0,\,1\}$ to be the characteristic function of $\M$, then 
\begin{align}\label{resnt}
    \Big|\sum_{mk = n\leqslant \sqrt T} \frac{r(m)\overline{r(n)}}{ k} (\log k)^{\ell}
\Big| \Big/ \Big(\sum_{n\leqslant \sqrt T }|r(n)|^2\Big) =  \frac{1}{|\M|}\sum_{\substack{ n \in \M\\ k|n}}\frac{(\log k)^{\ell}}{k} = \sum_{\substack{ k \in \M}}\frac{(\log k)^{\ell}}{k} \left( \frac{1}{|\M|}\sum_{\substack{ n \in \M\\ k|n}} 1 \right) \,.
\end{align}
Define $K:\, = \Big\{k\in \N\Big|\, k \leqslant \exp\left( \log_2 T \cdot \log_3 T\right)\,,~~P^{+}(k) \leqslant y\Big\}$ and its two subsets $K_1$ and $K_2$ to be
\begin{align*}
  &K_1:\, = \Big\{k\in K\Big|\, \sum_{i = 1}^w \alpha_i  \leqslant \frac{(\log_2 T)^3}{\log_3 T},  \,\, ~~\text{where} \,\,k\,\, \text{has the prime factorization as} ~~~  k = p_1 ^{\alpha_1} p_2 ^{\alpha_2} \cdots p_{w} ^{\alpha_{w}}\, \Big\}\,,  \\
  &K_2:\, = \Big\{k\in K\Big|\, \sum_{i = 1}^w \alpha_i  > \frac{(\log_2 T)^3}{\log_3 T},  \,\, ~~\text{where} \,\,k\,\, \text{has the prime factorization as} ~~~  k = p_1 ^{\alpha_1} p_2 ^{\alpha_2} \cdots p_{w} ^{\alpha_{w}}\, \Big\}\,.  
\end{align*}
Clearly, $K_1$ is a subset of $\M$. Let $k$ be any given element of $K_1$, then the inner sum in \eqref{resnt} tends to $1$, as $T \to \infty$. More precisely, we have
\begin{align}\label{closeTo1}
 \frac{1}{|\M|}\sum_{\substack{ n \in \M\\ k|n}} 1    \geqslant 1 - \frac{2}{\log_3 T}\,,\quad\quad \forall k \in K_1.
\end{align}
To see this, assume that $  k = p_1 ^{\alpha_1} p_2 ^{\alpha_2} \cdots p_{w} ^{\alpha_{w}}$. Then \[ \frac{1}{|\M|}\sum_{\substack{ n \in \M\\ k|n}} 1 = \frac{1}{b^w}\sum_{\substack{ n \in \M\\ k|n}} 1  = \frac{1}{b^w}\prod_{i = 1}^w \left(b - \alpha_i\right) = \prod_{i = 1}^w e^{\log \left( 1 - \frac{\alpha_i}{b}\right)}\,,\]
and \eqref{closeTo1} follows from the condition  $\sum_{i = 1}^w \alpha_i \leqslant   \frac{(\log_2 T)^3}{\log_3 T}\,.$

Now consider upper bounds for the  sum of reciprocals of elements of  $K_2$. By Rankin's trick and dropping conditions for $\alpha_i$, we have
\begin{align*}
    \sum_{k \in K_2} \frac{1}{k} \leqslant  \sum_{\alpha_1 = 0}^{\infty}   \sum_{\alpha_2 = 0}^{\infty} \cdots \sum_{\alpha_w = 0}^{\infty} \frac{1}{p_1 ^{\alpha_1} p_2 ^{\alpha_2} \cdots p_{w} ^{\alpha_{w}}} \left( \sum_{i = 1}^w  \alpha_i \right) \left(\frac{(\log_2 T)^3}{\log_3 T} \right)^{-1} \,.
\end{align*}
Next, we use the inequality $\sum_{i = 1}^w  \alpha_i  \leqslant \prod_{i = 1}^w (1 + \alpha_i)  $ and obtain
\begin{align*}
    \sum_{k \in K_2} \frac{1}{k} &\leqslant  \sum_{\alpha_1 = 0}^{\infty}  \frac{\alpha_1 + 1}{p_1^{\alpha_1}}  \sum_{\alpha_2 = 0}^{\infty}  \frac{\alpha_2 + 1}{p_2^{\alpha_2}}\cdots \sum_{\alpha_w = 0}^{\infty} \frac{\alpha_w + 1}{p_{w}^{\alpha_w}}  \left(\frac{(\log_2 T)^3}{\log_3 T} \right)^{-1}\\&\leqslant \frac{\log_3 T }{ (\log_2 T)^3  }  \prod_{i = 1}^w \left( \frac{1}{ 1 - \frac{1}{p_i}} \right)^2\\& \ll \frac{\log_3 T}{\log_2 T}\,,
\end{align*}
where in the last inequality we use the Mertens' theorem.  By the definition of $K_2$, when $k \in K_2$, we have $\log k \leqslant (\log_2 T) \cdot (\log_3 T)$. Thus we find that 
\begin{align}\label{K2}
 \sum_{k \in K_2} \frac{(\log k)^{\ell}}{k}  \ll  (\log_2 T)^{\ell -1} \cdot (\log_3 T)^{\ell + 1} \ll (\log_2 T)^{\ell} \cdot (\log_3 T) , \quad \forall  \ell \leqslant (\log_3 T) / (\log_4 T)\,.
\end{align}
In order to compute a lower bound for the outer sum in \eqref{resnt}, we first compute the sum over the set $K$, then by \eqref{K2} we restrict the sum to over its subset $K_1$, which is also a subset of $\M$.  Let $\RR = \exp\left( \log_2 T \cdot \log_3 T\right)\,.$ And we 
 keep in mind that $\ell \leqslant (\log_3 T) / (\log_4 T)$ in the following computations. 
 
 We split the sum into two parts as follows
\begin{align*}
 \sum_{k \in K} \frac{(\log k)^{\ell}}{k}   =   \sum_{k \leqslant y} \frac{(\log k)^{\ell}}{k} +  \sum_{ \substack{y < k \leqslant \RR\\P^{+}(k) \leqslant y}} \frac{(\log k)^{\ell}}{k}  = S_1 + S_2\,.
\end{align*}

The first sum is
\begin{align*}
  S_1 =\sum_{k \leqslant y} \frac{(\log k)^{\ell}}{k} = \left( \frac{1 }{\ell + 1} + o(1) \right) \left(\log_2 T \right)^{\ell+1}\,. 
\end{align*}
By partial summation, the second sum is 
\begin{align}\label{S2}
  S_2 = \frac{(\log R)^{\ell}}{R} \Psi (R, y) - \frac{(\log y)^{\ell}}{y} \Psi (y, y)  - \int_y^{\RR} \frac{d}{dx}   \left( \frac{(\log x)^{\ell}}{x}\right) \Psi(x, y) dx\,.
\end{align}
By \eqref{JianJin},  we have 
\begin{align}\label{Dickman}
   \Psi(x, y) = x\, \rho \left ( \frac{\log x}{\log y} \right) \left( 1 + O\left( \frac{\log_4 T}{\log_2 T} \right)\right)\,, \quad \text{uniformly for} \,\, y \leqslant x \leqslant R .
\end{align}
Applying \eqref{Dickman} into \eqref{S2}, and using \eqref{eqDickman} and \eqref{decay}, we obtain 
\begin{align*}
   S_2 =  \left ( \int_1^{\infty}  u^{\ell} \rho (u) du  + o(1) \right)  \left(\log_2 T \right)^{\ell+1}\,.
\end{align*}
We immediately get
\begin{align*}
 \sum_{k \in K} \frac{(\log k)^{\ell}}{k}   = S_1 + S_2  = \left ( \int_0^{\infty}  u^{\ell} \rho (u) du  + o(1) \right)  \left(\log_2 T \right)^{\ell+1}\,.
\end{align*}
Together with \eqref{K2}, we have
\begin{align}\label{K1}
 \sum_{k \in K_1} \frac{(\log k)^{\ell}}{k}   =   \left ( \int_0^{\infty}  u^{\ell} \rho (u) du  + o(1) \right)  \left(\log_2 T \right)^{\ell+1}\,.
\end{align}
Since $K_1$ is a subset of $\M$, we find that 
\begin{align*}
     \sum_{\substack{ k \in \M}}\frac{(\log k)^{\ell}}{k} \left( \frac{1}{|\M|}\sum_{\substack{ n \in \M\\ k|n}} 1 \right) \geqslant \sum_{\substack{ k \in K_1}}\frac{(\log k)^{\ell}}{k} \left( \frac{1}{|\M|}\sum_{\substack{ n \in \M\\ k|n}} 1 \right) \,.
\end{align*}
 By \eqref{resnt}, \eqref{closeTo1} and  \eqref{K1}, we are done.
\end{proof}
\section{Proof of theorem \ref{Main: sigma =  1}}
\begin{proof}

 By \cite[page 496]{DY}, we have
\begin{align*}
\max_{T\leqslant t\leqslant 2T}\left|\zeta^{(\ell)}\Big(1+it\Big)\right| & \geqslant \big(1 + O(T^{-1})\big)
  \Big|\sum_{mk = n\leqslant \sqrt T} \frac{r(m)\overline{r(n)}}{ k} (\log k)^{\ell}
\Big| \Big/ \Big(\sum_{n\leqslant \sqrt T }|r(n)|^2\Big) \\
&+\,O\Big( T^{-\frac{3}{2}}\, (\log T)^{\ell +1} \Big)  + O\Big(  (\log_2 T)^{\ell}\Big) \,. \end{align*}
    
By  Proposition \ref{maxRatio}, we finish the proof of Theorem \ref{Main: sigma =  1}.

\end{proof}

\section{Proof of theorem \ref{Lchi_lower} }

\begin{proof}

First note that we have the following approximation formula by  partial summation and  P\'{o}lya-Vinogradov inequality (\cite[Thm 9.18]{M2})
\begin{align}\label{approxD}
    L^{(\ell)}(1, \chi) = \sum_{k \leqslant N} \frac{\chi(k) (-\log k)^{\ell}}{k} + O\left( \frac{\sqrt{q} \log q (\log N)^{\ell}}{N} \right)\,, \quad \text{when} \quad \chi \neq \chi_0\,, \quad   \ell \leqslant \log N\,.
\end{align}
In order to use Soundararajan's resonance method \cite{So} to produce  extreme values, we define $V_2(q)$ and $V_1(q)$ as follows (also see \cite[page 129]{delaBT})
\begin{align*}
    V_2(q):\,  = \sum_{\chi \neq \chi_0}  (-1)^{\ell}  L^{(\ell)}(1, \chi; N) \left|R_{\chi}\right|^2\,, \quad
    V_1(q):\,  = \sum_{\chi\neq \chi_0}  \left|R_{\chi}\right|^2\,,
\end{align*}
where $L^{(\ell)}(1, \chi; N)$ and the resonator  $R_{\chi}$ are defined by 
\begin{align*}
      L^{(\ell)}(1, \chi; N):\, = \sum_{k \leqslant N} \frac{\chi(k) (-\log k)^{\ell}}{k} \,, \quad
    R_{\chi}:\,  = \sum_{a \leqslant A}\chi(a) r(a)\,.
\end{align*}
We chose $T = q^{\frac{1}{2}}$, $N = q^{\frac{3}{4}}$, $A = q^{\frac{1}{4}}$ and we let the function $r(n)$ to be defined  as in the proof of Proposition \ref{maxRatio}. By  orthogonality of characters, we have
\begin{align}\label{V1Shangjie}
      V_1(q) \leqslant \sum_{\chi}  \left|R_{\chi}\right|^2 \leqslant \phi(q) \sum_{a \leqslant A} r(a)\,. 
\end{align}
By Cauchy’s inequality, we have \[\left|R_{\chi_0}\right|^2 \leqslant A \sum_{a \leqslant A} r(a)\,.\]
Thus we can bound  $\left| L^{(\ell)}(1, \chi_0; N)\right|\cdot \left|R_{\chi_0}\right|^2$ by
\begin{align*}
    \leqslant (\log q)^{\ell + 1}  A \sum_{a \leqslant A} r(a).
\end{align*}
Above upper bound together with the orthogonality of characters gives that
\begin{align}\label{V2JianJin}
    V_2(q) = \phi(q) \sum_{mk = n \leqslant A}  \frac{(\log k)^{\ell} r(m)r(n)}{k}  +  O\left( (\log q)^{\ell + 1}    \right) \cdot  A \sum_{a \leqslant A} r(a) \,.
\end{align}
Combining \eqref{V2JianJin} with \eqref{V1Shangjie}, we have 
\begin{align}\label{ChiRatio}
    \max_{\chi \neq \chi_0} \left|L^{(\ell)}(1, \chi; N) \right| \geqslant \left| \frac{V_2(q)}{V_1(q)} \right| = \left(
    \sum_{mk = n \leqslant A}  \frac{(\log k)^{\ell} r(m)r(n)}{k}\right) \Big/\left( \sum_{a \leqslant A} r(a) \right) + O\left( \left(\log q\right)^{\ell+1} \right)\cdot q^{-\frac{3}{4} } \,.
\end{align}
By \eqref{ChiRatio}, \eqref{approxD} and  Proposition \ref{maxRatio}, we obtain \eqref{ExtremeL}.

\end{proof}


\section{Proof of Theorem \ref{RH_upperBound}}

\begin{proof}
Let $x_1 = \exp \left( (\log_2 T)^2 \right)$, $x_2 = T $, and  $y_j =  \log^2 T\log^{2} x_j\, (\log_2 T)^{12}$ for $j = 1,\, 2.$ Note that we have $\log y_1 \sim 2 \log_2 T,$ as $ T \to \infty.$
By taking  $\sigma=1$ and $\epsilon=(\log_2 T)^{-1}$ in \eqref{approx}, we have
$$
(-1)^{\ell}\zeta^{(\ell)}(1+i t)
= \sum_{k\leqslant T}\frac{(\log k)^{\ell}}{k^{1+i t}} + O\big((\log_2 T)^{\ell}\big),  \quad \forall t \in [2T,\, 5T]\,.
$$
We split the sum in the above approximation formula into two parts as follows
\begin{align*}
 \sum_{k \leqslant x_2} \frac{ (\log k)^{\ell}}{k^{1+it}} = \sum_{k \leqslant y_1 } \frac{ (\log k)^{\ell}}{k^{1+it}} + \sum_{y_1 < k  \leqslant x_2} \frac{ (\log k)^{\ell}}{k^{1+it}}.
\end{align*}
For the first sum, we have
\begin{align*}
\left| \sum_{k \leqslant y_1 } \frac{ (\log k)^{\ell}}{k^{1+it}} \right| \leqslant \sum_{k \leqslant y_1 } \frac{ (\log k)^{\ell}}{k} = \left( \frac{1}{\ell + 1}  + o(1) \right) \left(\log y_1 \right)^{\ell + 1} = \left( \frac{1 }{\ell + 1} + o(1) \right) \left(2\log_2 T \right)^{\ell+1}\,.
\end{align*}
For the second sum, by partial summation, we have
\begin{align*}
 \sum_{y_1 < k  \leqslant x_2} \frac{ (\log k)^{\ell}}{k^{1+it}} = \frac{ (\log x_2)^{\ell}}{x_2} \left(\sum_{k\leqslant x_2}  \frac{1}{k^{it}} \right)  -  \frac{ (\log y_1)^{\ell}}{y_1} \left(\sum_{k\leqslant y_1}  \frac{1}{k^{it}} \right)   +  \int_{y_1}^{x_2} \left(\sum_{n\leqslant x}  \frac{1}{n^{it}} \right) \frac{d}{dx}   \left( \frac{-(\log x)^{\ell}}{x}\right)  dx\,.  
\end{align*}
By \eqref{Upper_BySmooth} of Theorem \ref{RH_DY} and Lemma \ref{HTDickman}, we have
\begin{align*}
 \left| \frac{ (\log x_2)^{\ell}}{x_2} \left(\sum_{k\leqslant x_2}  \frac{1}{k^{it}} \right) \right| \ll \frac{(\log T)^{\ell}}{T} \Psi (T, \log^2 T (\log_2 T)^{20}) \ll (\log T)^{\ell}\, T^{-\frac{1}{2} + o(1)} \ll o(1) \cdot \left(\log_2 T \right)^{\ell+1}\,.
\end{align*}
Clearly, we have 
\begin{align*}
    \left|   \frac{ (\log y_1)^{\ell}}{y_1} \left(\sum_{k\leqslant y_1}  \frac{1}{k^{it}} \right) \right| \leqslant (\log y_1)^{\ell}\leqslant  o(1) \cdot \left(\log_2 T \right)^{\ell+1}\,.
\end{align*}
By \eqref{Approx_BySmooth} of Theorem \ref{RH_DY}, \eqref{eqDickman}  and Lemma \ref{HTDickman}, we have
\begin{align*}
 \left|   \int_{y_1}^{x_1} \left(\sum_{n\leqslant x}  \frac{1}{n^{it}} \right) \frac{d}{dx}   \left( \frac{-(\log x)^{\ell}}{x}\right)  dx \right| &\leqslant \int_{y_1}^{x_1} \left( 1 + o(1) \right) \Psi(x, y_1)  \frac{d}{dx}   \left( \frac{-(\log x)^{\ell}}{x}\right)  dx  \\
 & \leqslant \left ( \int_1^{\infty}  u^{\ell} \rho (u) du  + o(1) \right)  \left(2\log_2 T \right)^{\ell+1}\,.
\end{align*}
Again, by \eqref{Approx_BySmooth}, \eqref{eqDickman}  and Lemma \ref{HTDickman}, we have
\begin{align*}
 \left|   \int_{x_1}^{x_2} \left(\sum_{n\leqslant x}  \frac{1}{n^{it}} \right) \frac{d}{dx}   \left( \frac{-(\log x)^{\ell}}{x}\right)  dx \right| &\leqslant \int_{x_1}^{x_2} \left( 1 + o(1) \right) \Psi(x, y_2)  \frac{d}{dx}   \left( \frac{-(\log x)^{\ell}}{x}\right)  dx  \\
 & \leqslant  o(1) \cdot \left(\log_2 T \right)^{\ell+1}\,.
\end{align*}
As a result, we obtain
\begin{align*}
    \left| \zeta^{(\ell)}\left(1+it\right)\right| &\leqslant \left( \frac{1 }{\ell + 1} + o(1) \right) \left(2\log_2 T \right)^{\ell+1}  +  \left ( \int_1^{\infty}  u^{\ell} \rho (u) du  + o(1) \right)  \left(2\log_2 T \right)^{\ell+1}\\  & \leqslant \left( 2^{\ell + 1}\, \cl + o(1) \right) \left(\log_2 T \right)^{\ell + 1}\,.
\end{align*}
Since $t \in [2T,\, 5T]$, we are done.
\end{proof}

\section{Proof of Theorem \ref{RH_DY}}
The proof is almost the same as the proof of  the Granville-Soundararajan Theorem \cite[page 389--391]{largeGS} .  Only a few modifications are needed.
\begin{proof}

Define 
$$
\zeta(s; y) = \zeta(s) \prod_{p\leqslant y} \biggl(1-\frac{1}{p^s}\biggr)\,,
$$ so that $\zeta(s; y)$ is a meromorphic function on the whole plane, which is holomorphic everywhere except for a simple pole at $s = 1$.
When $|\text{Im}(s)|\leqslant T$    and  $ T + y + 3\leqslant t \leqslant T^{1000}$, we have $ y+ 3 \leqslant \text{Im}(s + it)\leqslant T + T^{1000}$. Note that
\begin{align*}
\log  \zeta(s + it ; y) &= \log  \zeta(s + it ) - \sum_{p \leqslant  y} \frac{1}{p^{s + it}} 
+ O\biggl(\sum_{p \leqslant y} \frac{1}{p^{2\text{Re}(s)}}\biggr),
\end{align*}
and so if $1 >$ Re$(s)\geqslant \frac{1}{2} +\frac{1}{\log y}$, \,$ T + y + 3\leqslant t \leqslant T^{1000}$,\, $ y\geqslant 2$  and 
$|\text{Im}(s)|\leqslant T$,
we get by Lemma \ref{approx_logZeta} 
$$\left|\log  \zeta(s + it ; y)\right| \leqslant C\log T \log^2 y, $$
where $C>0$ is some constant.  
Now suppose that $ x \in\mathbb{N} + \frac12$.
Let $u=\frac{\log x}{\log y}$ and put $c=1+\frac{1}{\log x}$.  By Perron's 
formula
\begin{align*}
\sum_{n\leqslant x}  \frac{1}{n^{it}} - \Psi(x,y;t) 
&=\frac{1}{2\pi i} 
\int_{c-i\infty}^{c+i\infty} 
\prod_{p\leqslant y}
\biggl(1-\frac{1}{p^{s+it}}\biggr)^{-1} 
\left( \exp \biggl( \log  \zeta(s + it ; y) \biggr)-1 \right) \frac {x^s}{s} ds \\ 
&=\sum_{j=1}^{[u]} \frac{1}{j!} \sum_{ \substack{n \leqslant  x/y^{j}  \\ P^{+}(n) \leqslant y   }  } \frac{1}{n^{it}}
\frac{1}{2\pi i} \int_{c-i\infty}^{c+i\infty} \biggl( \log  \zeta(s + it ; y) \biggr)^{j} 
\biggl( \frac{x}{n} \biggr)^s \frac{ds}{s}.
\end{align*}
Note that $\log \zeta(s ; y) = \sum_{n = 2}^{\infty} \Lambda_y(n) (\log n) ^{-1} n^{-s}$, where the generalized von Mangoldt function $\Lambda_y(\cdot)$ is defined as $\Lambda_y(n) = \log p$ if $n = p^{k}$ and $p > y$, otherwise $\Lambda_y(n) = 0.$ 

So we have $\left(\log  \zeta(s + it ; y)\right)^{j}/j! = \sum_{m=1}^{\infty} 
a_{j}(m,y)m^{-s - it}$, where $|a_j(m,y)|\leqslant 1$ for all $m$, $j$ and $y$.
All other
steps are the same as the proof  of the Granville-Soundararajan Theorem.
\end{proof}

\section{Proof of Theorem \ref{GRH_upperBound}}
\begin{proof}
Let $x_1 = \exp \left( (\log_2 q )^2 \right)$, $x_2 =  q^{\frac{3}{4}} $ , and  $ y_j = \log^2q\log^{2} x_j \, (\log_2 q )^{12}$ for $j = 1,\, 2.$  We will use the  approximation formula \eqref{approxD} for $L^{(\ell)}(1, \chi) $  and other steps are the same as the proof of Theorem  \ref{RH_upperBound}.
\end{proof}

\section{Proof of Theorem \ref{GSC_upper}}
\begin{proof}
Let $x_1 = \exp \left( (\log_2 q )^2 \right)$, $x_2 =  q^{\frac{3}{4}} $ , and  $ y_j = (\log q + \log^2 x_j) (\log_2 q )^A  $ for $j = 1,\, 2.$  We again use \eqref{approxD}  and other steps are the same as the proof of Theorem  \ref{RH_upperBound}. Note that now we have $\log y_1 \sim  \log_2 q ,$ as $ q \to \infty.$ Thus in the end, we obtain $\cl\left(\log_2 q\right)^{\ell+1}$ instead of $\cl\left(2\log_2 q\right)^{\ell+1}$ in Theorem  \ref{GRH_upperBound}\,.
\end{proof}

\section{Proof of Theorem \ref{DYC_Cor}}
\begin{proof}
For the upper bound, let $x_1 = \exp \left( (\log_2 T)^2 \right)$, $x_2 = T $, and  $y_j  = (\log T +\log^2 x_j) (\log_2 T)^A $ for $j = 1,\, 2.$ And other steps are the same as the proof of Theorem  \ref{RH_upperBound}. Combining with the lower bound, we are done.
\end{proof}

\section{A mixed conjecture}
Combining the Granville-Soundararajan Conjecture and Conjecture \ref{DYC}, we pose the following mixed conjecture.
\begin{conjecture}
There exists a constant $A>0$ such that
for any non-principal character $\chi$ (mod $q$), and for any
$1\leqslant x\leqslant \min\{q,\,T\}$, $2 T \leqslant t \leqslant 5T$, we have, uniformly,  
 
$$
\sum_{n\leqslant x} \frac{\chi(n)}{n^{it}} = \sum_{\substack{n \leqslant x \\ P^{+}(n) \leqslant y}}  \frac{\chi(n)}{n^{it}}  +o(\Psi(x,y;\chi_0)), \quad \text{as}\quad q \to \infty\,,T \to \infty\,,
$$
where $y= (\log qT +\log^2 x) (\log_2 qT)^A$.
\end{conjecture}

\section*{Acknowledgements}
I am grateful to Kannan Soundararajan   for several helpful discussions. I thank Shou-Wu Zhang
for pointing \cite{SC} to me. The work was carried out  when I was visiting the Stanford University. I thank the members of the institute for their hospitality. The work  was supported by the Austrian Science Fund (FWF), project W1230 and by the Simons  Grant of Kannan Soundararajan.

\end{document}